\documentclass[12pt]{article}
\usepackage{amsthm, amssymb, amsmath, latexsym}
\usepackage[utf8]{inputenc}
\usepackage[english]{babel}
\usepackage{longtable}

\begin{document}
\pagestyle{myheadings}
\thispagestyle{empty}
\setcounter{page}{1}

\newtheorem{definition}{Definition}
\newtheorem{proposition}{Proposition}
\newtheorem{theorem}{Theorem}
\newtheorem{lemma}{Lemma}
\newtheorem{corollary}{Corollary}
\newtheorem{remark}{Remark}
\theoremstyle{plain}
\mathsurround 2pt

\gdef\Aut{\mathop{\rm Aut}\nolimits}
\gdef\Ker{\mathop{\rm Ker}\nolimits}
\gdef\Im{\mathop{\rm Im}\nolimits}
\gdef\End{\mathop{\rm End}\nolimits}
\gdef\Inn{\mathop{\rm Inn}\nolimits}
\gdef\exp{\mathop{\rm exp}\nolimits}

\begin{center}
\textbf{\Large Groups of the nilpotency class $2$ of order $p^4$ as additive groups of local nearrings}
\end{center}

\begin{center}
\emph{Iryna Raievska, Maryna Raievska}\\
University of Warsaw, Warsaw, Poland;\\
Institute of Mathematics of National Academy of Sciences of Ukraine,\\
Kyiv, Ukraine\\
raeirina@imath.kiev.ua, raemarina@imath.kiev.ua
\end{center}

\date{}




\begin{abstract}
We consider groups of the nilpotency class $2$ of order $p^4$ which are the additive groups of local nearrings. It was shown that, for odd p, out of 6 of such groups 4 of them are the additive groups of local nearrings. Some examples of such nearrings are explicitly constructed.
\end{abstract}

\section{Introduction}

\

The question of finding groups that can be additive groups for the nearrings with identity is studied from the late 1960s. Some results in this direction were obtained in \cite{ClMal_66} and \cite{ClDoi_68}, where it was shown that the symmetric group $S_n$ with $n\geq 3$ and the alternating group $A_4$ cannot be an additive group of a nearring with identity, respectively. There is no nearring with identity whose additive group is isomorphic to the quaternion group $Q_8$~\cite{Clay_70}.

A study of local nearrings was first initiated in~\cite{MCJ_68} and it was found that the additive group of a finite zero-symmetric local nearring is a $p$-group. In~\cite{CM_68} it is shown that, up to an isomorphism, there exist $p-1$ local zero-symmetric nearrings with elementary abelian additive groups of order $p^2$, in which the subgroups of non-invertible elements have order $p$, that is, those nearrings which are not nearfields. Together with the fundamental paper~\cite{Z_36} and \cite{ClMal_66}, a complete description of all zero-symmetric local nearrings of order $p^2$ is obtained. The dihedral group $D_4$ of order $8$ cannot be the additive group of local nearrings~\cite{CM_71}. The existence of local nearrings on finite abelian $p$-groups is proved in~\cite{CM_70}, i.e. every non-cyclic abelian $p$-group of order $p^n>4$ is the additive group of a zero-symmetric local nearring which is not a ring. Also, it is established in \cite{IMYa_12} that an arbitrary non-metacyclic Miller--Moreno $p$-group of order $p^n>8$ is the additive group of some local nearring. Nearrings with identity and local nearrings on Miller--Moreno groups were studied in~\cite{IMYa_12}, \cite{IMYa_11} and \cite{IM_20}.

Boykett and N\"obauer~\cite{BN_98} classified all non-abelian groups of order less than $32$ that can be the additive groups of a nearring with identity and found the number of non-isomorphic nearrings with identity on such groups. The package SONATA~\cite{SONATA} of the computer algebra system~GAP~\cite{GAP} contains a library of all non-isomorphic nearrings of order at most $15$ and nearrings with identity of order up to $31$, among which $698$ are local.

However, the classification of nearrings of higher orders requires much more complex calculations. For local nearrings they were realized in the form of a new GAP package called LocalNR~\cite{LocalNR}. Its current version (not yet distributed with GAP) contains $37599$ local nearrings of order at most $361$, except orders $128$, $256$ and some of orders $32$, $64$ and $243$. We have already calculated some classes of local nearrings of orders $32$, $64$, $128$, $243$ and $625$.

However, it is not true that any finite group is the additive group of a nearring with identity. Therefore it is important to determine such groups and to classify some classes of nearrings with identity on these groups, for example, local nearrings.

In \cite{IM_2022} it was shown that on each group of order $p^3$ with $p>2$ there exists a local nearring. Moreover, lower bounds for the number of local nearrings on groups of order $p^3$ are obtained. It is established that on each non-metacyclic non-abelian or metacyclic abelian groups of order $p^3$ there exist at least $p+1$ non-isomorphic local nearrings. In~ \cite{IM_2021} it is proved that, up to an isomorphism, there exist at least $p$ local nearrings on elementary abelian additive groups of order $p^3$, which are not nearfields.

The next natural step is to investigate groups of order $p^4$ as the additive groups of local nearrings. In this paper we consider groups of the nilpotency class $2$ of order $p^4$ which are the additive group of local nearrings.

\section{Preliminaries}

\

We will give the basic definitions.

\begin{definition}
A non-empty set $R$ with two binary operations $``+"$ and $``\cdot"$ is a \textbf{nearring} if:
\begin{description}
  \item[1)] $(R,+)$ is a group with neutral element $0${\rm ;}
  \item[2)] $(R,\cdot)$ is a semigroup{\rm ;}
  \item[3)] $x\cdot (y+z)=x\cdot y+x\cdot z$ for all $x$, $y$, $z\in R$.
\end{description}
Such a nearring is called a left nearring. If axiom 3) is replaced by an axiom $(x+y)\cdot z = x\cdot z + y\cdot z$ for all $x$, $y$, $z\in R$, then we get a right
nearring.
\label{nr}
\end{definition}

The group $(R,+)$ of a nearring $R$ is denoted by $R^+$ and called the {\em additive group} of $R$. It is easy to see that for each subgroup $M$ of $R^+$ and for each element $x\in R$ the set $xM=\{x\cdot y|y\in M\}$ is a subgroup of $R^+$ and in particular $x\cdot 0=0$. If in addition $0\cdot x=0$ for all $x\in R$, then the nearring $R$ is called {\em zero-symmetric}. Furthermore, $R$ is a {\em nearring with an identity} $i$  if the semigroup $(R,\cdot)$ is a monoid with identity element $i$. In the latter case the group of all invertible elements of the monoid $(R,\cdot)$ is denoted by $R^*$ and called the {\em multiplicative group} of $R$.  A subgroup $M$ of $R^+$ is called $R^*$-{\em invariant}, if $rM\leq M$ for each $r\in R^*$, and $(R,R)$-{\em subgroup}, if $xMy\subseteq M$ for arbitrary $x$, $y\in R$.\medskip

The following assertion is well-known (see, for instance, \cite{ClMal_66}, Theorem 3).

\begin{lemma}\label{exponent}
The exponent of the additive group of a finite nearring $R$ with identity $i$ is equal to the additive order of $i$ which coincides with the additive order of every invertible element of $R$.

\end{lemma}

\begin{definition}
A nearring $R$ with identity is called \textbf{local} if the set $L$ of all non-invertible elements of $R$ forms a subgroup of the additive group $R^{+}$ and a \textbf{nearfield}, if $L=0$.
\end{definition}

Through this paper $L$ will denote the subgroup of non-invertible elements of $R$.

The following lemma characterizes the main properties of finite local nearrings (see~\cite{AHS_2004}, Lemma~3.2).

\begin{lemma}\label{prop}
Let $R$ be a local nearring with identity $i$. Then the following statements hold{\rm :}
\begin{description}
\item[1)] $L$ is an $(R,R)$-subgroup of $R^{+}${\rm ;}
\item[2)] each proper $R^*$-invariant subgroup of $R^+$ is contained in $L${\rm ;}
\item[3)] the set $i+L$ forms a subgroup of the multiplicative group $R^*$.
\end{description}
\end{lemma}

Finite local nearrings with a cyclic subgroup of non-invertible elements are described in~\cite[Theorem~1]{RIM_2}.

\begin{theorem}\label{theorem_2}
Let $R$ be a local nearring of order $p^n$ with {$n\!> 1$} whose subgroup $L$ is cyclic and non-trivial. Then the additive group $R^+$ is either cyclic or is an elementary abelian group of order $p^2$. In the first case, $R$ is a commutative local ring, which is isomorphic to residual ring $\mathbb Z/p^n\mathbb Z$ with $n\ge 2$, in the other case there exist $p$ non-isomorphic such nearrings $R$ with $|L|=p$, from which $p-1$ are zero-symmetric nearrings and their multiplicative groups $R^{*}$ are isomorphic to a semidirect product of two cyclic subgroups of orders $p$ and $p-1$.
\end{theorem}

The following theorem was proved by Maxson in \cite{CM_68} (Theorem 2.1).

\begin{theorem}\label{theorem_3}
If $R$ is a finite local nearring which is not a nearfield, then $|R|<|L|^2$.
\end{theorem}

As a consequence of Theorems~\ref{theorem_2} and \ref{theorem_3} we have the following result.

\begin{corollary}
Let $R$ be a local nearring of order $p^4$ with non-abelian additive group and is not a nearfield. Then the subgroup of non-invertible elements $L$ is a non-cyclic group of order $p^3$ or $p^2$.
\label{cor_1}
\end{corollary}

We recall the following definition.

\begin{definition}
A finite non-abelian group whose proper subgroups are abelian is called a \textbf{Miller--Moreno group}.
\end{definition}

\section{Groups of the nilpotency class $2$ of order~$p^4$}

\

We will consider groups of the nilpotency class $2$ of order $p^4$.

Let $[n,i]$ be the $i$-th group of order $n$ in the SmallGroups library in the computer system algebra GAP. We denote by $C_n$ the cyclic group of order $n$.

It is an easy exercise for example in GAP to get the following assertion.

\begin{remark}\label{2groups}
There are $6$ groups of the nilpotency class $2$ of order $2^4=16$, which are:
\begin{enumerate}
  \item $(C_4 \times C_2) \rtimes C_2$ $[ 16, 3 ]${\rm ;}
  \item $C_4 \rtimes C_4$ $[ 16, 4 ]${\rm ;}
  \item $C_8 \rtimes C_2$ $[ 16, 6 ]${\rm ;}
  \item $C_2 \times D_8$ $[ 16, 11 ]${\rm ;}
  \item $C_2 \times Q_8$ $[ 16, 12 ]${\rm ;}
  \item $(C_4 \times C_2) \rtimes C_2$ $[ 16, 13 ]$.
\end{enumerate}
\end{remark}

The following theorem contains the classification of groups of the nilpotency class $2$ of order $p^4$, where $p$ is an odd prime (see, \cite{Burnside_1897} and, for example, \cite{Al-Hasanat_Almazaydeh_2022}).

\begin{theorem}\label{groups}
There are $6$ groups of the nilpotency class $2$ of order $p^4$, where $p$ is an odd prime, which are:
\begin{itemize}
  \item $G_1 = \langle a, b \colon a^{p^2}=b^p=[a, [a, b]] = [b, [a, b]] = [a, b]^p = e\rangle = (C_{p^2}\times C_p)\rtimes C_p${\rm ;}
  \item $G_2 = \langle a, b \colon a^{p^2}=b^{p^2}= e, [b, a] = b^p\rangle= C_{p^2}\rtimes C_{p^2}${\rm ;}
  \item $G_3 = \langle a, b \colon a^{p^3}= b^p = e, [b, a] = a^{p^2}\rangle= C_{p^3}\rtimes C_p${\rm ;}
  \item $G_4 = \langle a, b, d \colon a^p = b^p = c^p = d^p = [a, c] = [b, c] = [a, d] = [b, d] = [c, d] = e\rangle = C_p \times ((C_p \times C_p)\rtimes C_p)$, where $c = [a, b]${\rm ;}
  \item $G_5 = \langle a, b, c \colon a^{p^2}= b^p = c^p = [a, c] = [b, c] = e, [b, a] = a^p\rangle = C_p\times (C_{p^2}\rtimes C_p)${\rm ;}
  \item $G_6 = \langle a, b, c \colon a^{p^2}= b^p = c^p = [a, b] = [a, c] = e, [c, b] = a^p\rangle = (C_{p^2}\times C_p)\rtimes  C_p$.
\end{itemize}
\end{theorem}

\section{Groups of the nilpotency class $2$ of order $16$ and local nearrings}

\

As was mentioned above a library of all non-isomorphic nearrings with identity of order up to $31$ are contained in the package SONATA, and so all non-isomorphic local nearrings of order $16$ (see~\cite{BN_98}).

\begin{lemma}\label{2addgroups}
The following groups of the nilpotency class $2$ and only they are the additive groups of local nearrings of order $16$\rm{:}
\begin{enumerate}
  \item $(C_4 \times C_2) \rtimes C_2$ $[ 16, 3 ]${\rm ;}
  \item $C_4 \rtimes C_4$ $[ 16, 4 ]${\rm ;}
  \item $C_8 \rtimes C_2$ $[ 16, 6 ]${\rm ;}
  \item $C_2 \times Q_8$ $[ 16, 12 ]$.
\end{enumerate}
\end{lemma}

Let $n(G)$ be the number of all non-isomorphic local nearrings $R$ whose additive group $R^+$ is isomorphic to the group~$G$.

\begin{center}
\begin{tabular}{|c|c|}
\hline
$StructureDescription(R^+)$ & $n(R^+)$\\
 \hline
$(C_4 \times C_2) \rtimes C_2$        & 37 \\
 \hline
$C_4 \rtimes C_4$                     & 24 \\
 \hline
$C_8 \rtimes C_2$                     & 33 \\
 \hline
$C_2 \times Q_8$                      & 2 \\
 \hline
\end{tabular}
\end{center}

\subsection{The groups $G_1$, $G_2$ and $G_3$ and local nearrings}

\

The groups $G_1$, $G_2$ and $G_3$ from Theorem~\ref{groups} are Miller--Moreno groups.

Due to \cite{IYa_12} $G_2$ and $G_3$ are the groups $G(p^2,p^2)$ and $G(p^3,p)$, respectively (see, for example, Lemma~2~\cite{IYa_12}). Therefore, by Theorem 2~\cite{IYa_12} there exists a local nearring $R$ whose additive group $R^+$ is isomorphic to $G_3$. As a consequence, there does not exit a local nearring on the additive group $G_2$.

Let $R$ be a local nearring whose additive group of $R^+$ is isomorphic to $G_2$. Then $R^{+}=\langle a\rangle +\langle b \rangle$ for some elements $a$ and $b$ of $R$ satisfying the relations $a{p}^3=0$, $bp=0$ and $-b+a+b=a(1-p^2)$. In particular, each element $x\in R$ is uniquely written in the form $x=ax_1+bx_2$ with coefficients $0\le x_1<p^3$ and $0\le x_2<p$.

The formula for multiplying elements of local nearrings on Miller--Moreno metacyclic groups is defined in~\cite{IYa_12}. The multiplication formula for arbitrary elements of a zero-symmetric local nearring  on $G_2$ is given in the proving of~\cite[Theorem~2]{IYa_12}, namely:
$$x\cdot y=a(x_1y_1+p^2x_1x_2\binom{y_1}{2})+b(x_2y_1+\beta(x)y_2),$$
where $\beta(x)=\left\{
               \begin{array}{ll}
                 1, & if~\hbox{$x_1\not\equiv 0\; (\!\!\mod p\;)$;}\\
                 0, & if~\hbox{$x_1\equiv 0\; (\!\!\mod p\;)$.}
               \end{array}
             \right.$

\textbf{Example 1.} Let $G\cong C_{27} \rtimes C_3$. If $x=ax_1+bx_2$ and $y=ay_1+by_2\in G$ and $(G,+, \cdot)$ is a local nearring, then as above $``\cdot"$ can be the following multiplication:

$$x\cdot y=a(x_1y_1+9x_1x_2\binom{y_1}{2})+b(x_2y_1+\beta(x)y_2),$$
where $\beta(x)=\left\{
               \begin{array}{ll}
                 1, & if~\hbox{$x_1\not\equiv 0\; (\!\!\mod 3\;)$;}\\
                 0, & if~\hbox{$x_1\equiv 0\; (\!\!\mod 3\;)$.}
               \end{array}
             \right.$

A computer program verified that the nearring obtained in Example~1 is indeed a local nearring, is deposited on GitHub:

\verb+https://github.com/raemarina/Examples/blob/main/LNR_81-6.txt+

From the package LocalNR and \cite{zenode_625} we have the following number of all non-isomorphic zero-symmetric local nearrings on $G_3$ of orders $81$ and $625$.

\begin{center}
\begin{tabular}{|c|c|}
\hline
$StructureDescription(R^+)$ & $n(R^+)$\\
 \hline
$C_{27} \rtimes C_3$                     & 10 \\
 \hline
$C_{125} \rtimes C_5$                    & 5  \\
 \hline
\end{tabular}
\end{center}

Analogously, $G_1$ is the group $G(p^2,p,p)$ according to \cite{IMYa_12}. Hence, by Theorem~3~\cite{IMYa_12} there exist a local nearring whose additive group is isomorphic to $G_1$. Since $G_1$ is a Miller--Moreno non-metacyclic group, using~\cite[Theorem~3]{IMYa_12}, for arbitrary elements $x=ax_1+bx_2+cx_3$ and $y=ay_1+by_2+cy_3$ of $G_1$ we obtain the following multiplication formula:
$$x\cdot y=a(x_1y_1+p^kx_2y_2)+b(x_2y_1+x_1y_2)+c(-x_1x_2\binom{y_1}{2}+x_3y_1+x_1^2y_3),$$
where $k=1,~2$.

It is easy to see that $R = (G_1,+,\cdot)$ is a non-zero-symmetric local nearring.

\textbf{Example 2.} Let $G\cong (C_9\times C_3)\rtimes C_3$. If $x=ax_1+bx_2+cx_3$ and $y=ay_1+by_2+cy_3\in G$ and $(G,+, \cdot)$ is a local nearring, then as above $``\cdot"$ can be one of the following multiplications:
\begin{description}
  \item[(1)] $x\cdot y=a(x_1y_1+3x_2y_2)+b(x_2y_1+x_1y_2)+c(-x_1x_2\binom{y_1}{2}+x_3y_1+x_1^2y_3)${\rm ;}
  \item[(2)] $x\cdot y=a(x_1y_1)+b(x_2y_1+x_1y_2)+c(-x_1x_2\binom{y_1}{2}+x_3y_1+x_1^2y_3)$.
\end{description}

A computer program verified that the nearring obtained in Example 2 is indeed a local nearring, is deposited on GitHub:

\verb+https://github.com/raemarina/Examples/blob/main/LNR_81-3.txt+

From the package LocalNR and \cite{zenode_625} we have the following number of all non-isomorphic zero-symmetric local nearrings on $G_1$ of orders $81$ and $625$.

\begin{center}
\begin{tabular}{|c|c|}
\hline
$StructureDescription(R^+)$ & $n(R^+)$\\
 \hline
$(C_9\times C_3)\rtimes C_3$                     & 46 \\
 \hline
$(C_{25}\times C_5)\rtimes C_5$                  & 154 \\
 \hline
\end{tabular}
\end{center}

\subsection{The group $G_4$}

\

Let $G_4$ be additively written group from Theorem~\ref{groups}. Then $G_4=\langle a\rangle +\langle b \rangle +\langle c \rangle+\langle d \rangle$ for some elements $a$, $b$, $c$ and $d$ of $R$ satisfying the relations $ap=0$, $bp=0$, $cp=0$, $dp=0$, $a+b=b+a+c$, $a+c=c+a$, $b+c=c+b$, $a+d=d+a$, $b+d=d+b$ and $c+d=d+c$.

\begin{lemma}\label{G4.1}
For arbitrary integers $k$ and $l$ in the group $G_4$ the equalities $-ak-bl+ak+bl=c(kl)$ and $bl+ak=-c(kl)+ak+bl$ hold.
\end{lemma}

\begin{proof}
Since $-a-b+a+b=c$, we get $-b+a+b=a+c$. Then
$$-bl+ak+bl=(a+cl)k=ak+c(kl).$$
Therefore, $-ak-bl+ak+bl=c(kl)$ and, so $bl+ak=-c(kl)+ak+bl$.
\end{proof}

\begin{lemma}\label{G4.2}
For any natural numbers $k$, $l$, $n$, $m$ and $r$ in the group $G_4$ the equality $(ak+bl+cm+dn)r=a(kr)+b(lr)+c(mr-kl\binom{r}{2})+d(nr)$ holds.
\end{lemma}

\begin{proof}
The proof will be carried out by induction on $r$. For $r=1$ the equality is valid.
Let for $r$ the equality hold, i.e.
$$(ak+bl+cm+dn)r=a(kr)+b(lr)+c(mr-kl\binom{r}{2})+d(nr).$$
Let us prove the equality for $r+1$:
$$\begin{array}{l}
(ak+bl+cm+dn)(r+1)=\\
\qquad{}=a(kr)+b(lr)+ak+bl+c(kl\binom{r}{2})+c(m(r+1))+d(n(r+1))=\\
\qquad{}=a(k(r+1))+b(l(r+1))+c(-klr)+c(-kl\binom{r}{2})+\\
\qquad{}+c(m(r+1))+d(n(r+1))=\\
\qquad{}=a(k(r+1))+b(l(r+1))+c(m(r+1)-kl(r+\binom{r}{2})+d(n(r+1))=\\
\qquad{}=a(k(r+1))+b(l(r+1))+c(m(r+1)-kl\binom{r+1}{2})+d{n(r+1)}.
\end{array}$$
Therefore, the equality is valid for any $r$.
\end{proof}

\subsection{Nearrings with identity whose\\
additive groups are isomorphic to $G_4$}

\

Let $R$ be a nearring with identity whose additive group $R^+$ is isomorphic to $G_4$. Then $R^{+}=\langle a\rangle +\langle b \rangle +\langle c \rangle+\langle d \rangle$ for some elements $a$, $b$, $c$ and $d$ of $R$ satisfying the relations $ap=0$, $bp=0$, $cp=0$, $dp=0$, $a+b=b+a+c$, $a+c=c+a$, $b+c=c+b$, $a+d=d+a$, $b+d=d+b$ and $c+d=d+c$. In particular, each element $x\in R$ is uniquely written in the form $x=ax_1+bx_2+cx_3+dx_4$ with coefficients $0\le x_1<p$, $0\le x_2<p$, $0\le x_3<p$ and $0\le x_4<p$.

Since the order of the element $a$ is equal to the exponent of group $G$, then by Lemma~\ref{exponent} we can assume that $a$ is an identity of $R$, i.e. $ax=xa=x$ for each $x\in R$. Furthermore, for each $x\in R$ there exist coefficients $\alpha(x)$, $\beta(x)$, $\gamma(x)$, $\varphi(x)$, $\lambda(x)$, $\mu(x)$, $\nu(x)$ and $\psi(x)$ such that $xb=a\alpha(x)+b\beta(x)+c\gamma(x)+d\varphi(x)$ and $xd=a\lambda(x)+b\mu(x)+c\nu(x)+d\psi(x)$. It is clear that they are uniquely defined modulo $p$, so that some mappings $\alpha\colon R\to \mathbb Z_{p}$, ${\beta: R\rightarrow Z_{p}}$, ${\gamma: R\rightarrow Z_p}$, ${\varphi: R\rightarrow Z_p}$, ${\lambda: R\rightarrow Z_p}$, ${\mu: R\rightarrow Z_p}$, ${\nu: R\rightarrow Z_p}$ and ${\psi: R\rightarrow Z_p}$ are determined.

\begin{lemma}\label{G4.3}
Let $R$ be a nearring with identity whose additive group $R^+$ is isomorphic to $G_4$. If $a$ coincides with identity element of $R$, $x=ax_1+bx_2+cx_3+dx_4$, $y=ay_1+by_2+cy_3+dy_4\in R$, $xb=a\alpha(x)+b\beta(x)+c\gamma(x)+d\varphi(x)$ and $xd=a\lambda(x)+b\mu(x)+c\nu(x)+d\psi(x)$, then
$$\begin{array}{c}
xy=a(x_1y_1+\alpha(x)y_2+\lambda(x)y_4)+b(x_2y_1+\beta(x)y_2+\mu(x)y_4)+\\
+c(-x_1x_2\binom{y_1}{2}-\alpha(x)\beta(x)\binom{y_2}{2}-x_2\alpha(x)y_1y_2-\lambda(x)\mu(x)\binom{y_4}{2}-x_2\alpha(x)y_3+\\
+x_3y_1+\gamma(x)y_2+x_1\beta(x)y_3+\nu(x)y_4)+d(x_4y_1+\varphi(x)y_2+\psi(x)y_4).\end{array}$$
Moreover, for the mappings\\
$\begin{array}{c}
\alpha\colon R\to \mathbb Z_{p}, {\beta: R\rightarrow Z_{p}}, {\gamma: R\rightarrow Z_p}, {\varphi: R\rightarrow Z_p},{\lambda: R\rightarrow Z_p}, {\mu: R\rightarrow Z_p},\\
{\nu: R\rightarrow Z_p}~\mbox{and}~{\psi: R\rightarrow Z_p}~\mbox{the following statements hold}{\rm :}\end{array}$
\begin{description}
  \item[\rm{(0)}] $\alpha(0)\equiv 0\; (\!\!\mod p)$, $\beta(0)\equiv 0\; (\!\!\mod p)$, ${\gamma(0)\equiv 0\; (\!\!\mod p)}${\rm ;}\\
   ${\varphi(0)\equiv 0\; (\!\!\mod p)}$, ${\lambda(0)\equiv 0\; (\!\!\mod p)}$, ${\mu(0)\equiv 0\; (\!\!\mod p)}$, \\
    ${\nu(0)\equiv 0\; (\!\!\mod p)}$ and ${\psi(0)\equiv 0\; (\!\!\mod p)}$ if and only if the nearring $R$ is zero-symmetric{\rm ;}
  \item[\rm{(1)}] $\alpha(xy)\equiv x_1\alpha(y)+\alpha(x)\beta(y)+\lambda(x)\varphi(y)\; (\!\!\mod p\;){\rm ;}$
  \item[\rm{(2)}] $\beta(xy)\equiv x_2\alpha(y)+\beta(x)\beta(y)+\mu(x)\varphi(y)\; (\!\!\mod p\;){\rm ;}$
  \item[\rm{(3)}] $\gamma(xy)\equiv -x_1x_2\binom{\alpha(y)}{2}-\alpha(x)\beta(x)\binom{\beta(y)}{2}-x_2\alpha(x)\alpha(y)\beta(y)-\lambda(x)\mu(x)\binom{\varphi(y)}{2}-$\\
      $-x_2\alpha(x)\gamma(y)+x_3\alpha(y)+\gamma(x)\beta(y)+x_1\beta(x)\gamma(y)+\nu(x)\varphi(y)\; (\!\!\mod p\;){\rm ;}$
  \item[\rm{(4)}] $\varphi(xy)\equiv x_4\alpha(y)+\varphi(x)\beta(y)+\psi(x)\varphi(y)\; (\!\!\mod p\;){\rm ;}$
  \item[\rm{(5)}] $\lambda(xy)\equiv x_1\lambda(y)+\alpha(x)\mu(y)+\lambda(x)\psi(y)\; (\!\!\mod p\;){\rm ;}$
  \item[\rm{(6)}] $\mu(xy)\equiv x_2\lambda(y)+\beta(x)\mu(y)+\mu(x)\psi(y)\; (\!\!\mod p\;){\rm ;}$
  \item[\rm{(7)}] $\nu(xy)\equiv -x_1x_2\binom{\lambda(y)}{2}-\alpha(x)\beta(x)\binom{\mu(y)}{2}-x_2\alpha(x)\lambda(y)\mu(y)-\lambda(x)\mu(x)\binom{\psi(y)}{2}-$\\
      $-x_2\alpha(x)\nu(y)+x_3\lambda(y)+\gamma(x)\mu(y)+x_1\beta(x)\nu(y)+\nu(x)\psi(y)\; (\!\!\mod p\;){\rm ;}$
  \item[\rm{(8)}] $\psi(xy)\equiv x_4\lambda(y)+\varphi(x)\mu(y)+\psi(x)\psi(y)\; (\!\!\mod p\;)$.
\end{description}
\end{lemma}

\begin{proof}
Since $0\cdot a=a\cdot 0=0$, it follows that $R$ is a zero-symmetric nearring if and only if
$$0=0\cdot b=a\alpha(0)+b\beta(0)+c\gamma(0)+d\varphi(0)$$
and
$$0=0\cdot d=a\lambda(0)+b\mu(0)+c\nu(0)+d\psi(0).$$
Equivalently we have

$\begin{array}{l} \alpha(0)\equiv 0\; (\!\!\mod p),~{\beta(0)\equiv 0\; (\!\!\mod p)},~{\gamma(0)\equiv 0\; (\!\!\mod p)},~{\varphi(0)\equiv 0\; (\!\!\mod p)},\\
~{\lambda(0)\equiv 0\; (\!\!\mod p)},~{\mu(0)\equiv 0\; (\!\!\mod p)},~{\nu(0)\equiv 0\; (\!\!\mod p)}~\mbox{and}~{\psi(0)\equiv 0\; (\!\!\mod p)}.\end{array}$

Moreover, since $c=-a-b+a+b$ and the left distributive law we have $0\cdot c=-0\cdot a-0\cdot b+0\cdot a+0\cdot b=0$, whence $$0\cdot x=0\cdot (ax_1+bx_2+cx_3+dx_4)=(0\cdot a)x_1+(0\cdot b)x_2+(0\cdot c)x_3+(0\cdot d)x_4=0.$$ So that statement (0) holds.

Further, using Lemma~\ref{G4.1}, we derive
$$\begin{array}{c}
xc=-xa-xb+xa+xb=-cx_3-bx_2-ax_1-c\gamma(x)-\\
-b\beta(x)-a\alpha(x)+ax_1+bx_2+cx_3+a\alpha(x)+b\beta(x)+c\gamma(x)=\\
=-bx_2-ax_1-b\beta(x)-a\alpha(x)+ax_1+bx_2+a\alpha(x)+b\beta(x)=\\
=-bx_2+cx_1\beta(x)-b\beta(x)-ax_1-a(\alpha(x)-x_1)+bx_2+a\alpha(x)+b\beta(x)=\\
=cx_1\beta(x)-b(x_2+\beta(x))-a\alpha(x)+bx_2+a\alpha(x)+b\beta(x)=\\
=cx_1\beta(x)-b(x_2+\beta(x))-a\alpha(x)-cx_2\alpha(x)+a\alpha(x)+bx_2+b\beta(x)=\\
=c(x_1\beta(x)-x_2\alpha(x))-b(x_2+\beta(x))+bx_2+b\beta(x)=c(x_1\beta(x)-x_2\alpha(x)).
\end{array}$$

Further, using the left distributive law, we obtain
$$\begin{array}{c}
xy=(ax_1+bx_2+cx_3+dx_4)y_1+(a\alpha(x)+b\beta(x)+c\gamma(x)+d\varphi(x))y_2+\\
+(c(x_1\beta(x)-x_2\alpha(x)))y_3+(a\lambda(x)+b\mu(x)+c\nu(x)+d\psi(x))y_4.
\end{array}$$
By Lemma~\ref{G4.2}, we get
$$\begin{array}{c}
(ax_1+bx_2+cx_3+dx_4)y_1=ax_1y_1+bx_2y_1+c(x_3y_1-x_1x_2\binom{y_1}{2})+dx_4y_1,\\
(a\alpha(x)+b\beta(x)+c\gamma(x)+d\varphi(x))y_2=a\alpha(x)y_2+b\beta(x)y_2+\\
+c(\gamma(x)y_2-\alpha(x)\beta(x)\binom{y_2}{2})+d\varphi(x)y_2,\\
(a\lambda(x)+b\mu(x)+c\nu(x)+d\psi(x))y_4=a\lambda(x)y_4+b\mu(x)y_4+\\
+c(\nu(x)y_4-\lambda(x)\mu(x)\binom{y_4}{2})+d\psi(x)y_4.\\
\end{array}$$
By Lemma~\ref{G4.2}, we have
$$bx_2y_1+a\alpha(x)y_2=a\alpha(x)y_2+bx_2y_1-cx_2\alpha(x)y_1y_2,$$
and
$$b\beta(x)y_2+a\lambda(x)y_4=a\lambda(x)y_4+b\beta(x)y_2-c\lambda(x)\beta(x)y_2y_4.$$
Hence and using the left distributive law, we have
$$\begin{array}{c}
xy=a(x_1y_1+\alpha(x)y_2+\lambda(x)y_4)+b(x_2y_1+\beta(x)y_2+\mu(x)y_4)+c(-x_1x_2\binom{y_1}{2}-\\
-\alpha(x)\beta(x)\binom{y_2}{2}-x_2\alpha(x)y_1y_2-\lambda(x)\mu(x)\binom{y_4}{2}-x_2\alpha(x)y_3+\\
+x_3y_1+\gamma(x)y_2+x_1\beta(x)y_3+\nu(x)y_4)+d(x_4y_1+\varphi(x)y_2+\psi(x)y_4).
\end{array}$$
The associativity of multiplication in $R$ implies that for all $x$, $y\in R$
$$(xy)b=x(yb)\leqno 1)$$ and
$$(xy)d=x(yd).\leqno 2)$$
According to $xb=a\alpha(x)+b\beta(x)+c\gamma(x)+d\varphi(x)$, we obtain $$(xy)b=a\alpha(xy)+b\beta(xy)+c\gamma(xy)+d\varphi(xy)\leqno 3)$$ and
$yb=a\alpha(y)+b\beta(y)+c\gamma(y)+d\varphi(y)$. Substituting the last equation to the right part of equality 1), we also have
$$\begin{array}{c}
x(yb)=a(x_1\alpha(y)+\alpha(x)\beta(y)+\lambda(x)\varphi(y))+b(x_2\alpha(y)+\beta(x)\beta(y)+\mu(x)\varphi(y))+\\
+c(-x_1x_2\binom{\alpha(y)}{2}-\alpha(x)\beta(x)\binom{\beta(y)}{2}-x_2\alpha(x)\alpha(y)\beta(y)-\\
-\lambda(x)\mu(x)\binom{\varphi(y)}{2}-x_2\alpha(x)\gamma(y)+x_3\alpha(y)+\gamma(x)\beta(y)+x_1\beta(x)\gamma(y)+\\
+\nu(x)\varphi(y))+d(x_4\alpha(y)+\varphi(x)\beta(y)+\psi(x)\varphi(y)).
\end{array}\leqno 4)$$
Since equality 1) implies the congruence of the corresponding coefficients in formulas 3) and 4), we obtain statements (1)--(4).

Next, according to $xd=a\lambda(x)+b\mu(x)+c\nu(x)+d\psi(x)$ instead of $y$ in equality~2), we get $$(xy)d=a\lambda(xy)+b\mu(xy)+c\nu(xy)+d\psi(xy)\leqno 5)$$ and
$yd=a\lambda(y)+b\mu(y)+c\nu(y)+d\psi(y)$. Substituting the last equation to the right part of equality 2), we also have
$$\begin{array}{c}
x(yd)=a(x_1\lambda(y)+\alpha(x)\mu(y)+\lambda(x)\psi(y))+b(x_2\lambda(y)+\beta(x)\mu(y)+\mu(x)\psi(y))+\\
+c(-x_1x_2\binom{\lambda(y)}{2}-\alpha(x)\beta(x)\binom{\mu(y)}{2}-x_2\alpha(x)\lambda(y)\mu(y)-\\
-\lambda(x)\mu(x)\binom{\psi(y)}{2}-x_2\alpha(x)\nu(y)+x_3\lambda(y)+\gamma(x)\mu(y)+x_1\beta(x)\nu(y)+\\
+\nu(x)\psi(y))+d(x_4\lambda(y)+\varphi(x)\mu(y)+\psi(x)\psi(y)).
\end{array}\leqno 6)$$
Finally, comparing the coefficients under $a$, $b$, $c$ and $d$ in formulas 5) and 6), we derive statements~(5)--(8) of the lemma.
\end{proof}

\subsection{Local nearrings whose additive groups\\
are isomorphic to $G_4$}

\

Let $R$ be a local nearring whose additive group $R^+$ is isomorphic to $G_4$. Then $R^{+}=\langle a\rangle +\langle b \rangle +\langle c \rangle+\langle d \rangle$ for some elements $a$, $b$, $c$ and $d$ of $R$ satisfying the relations $ap=0$, $bp=0$, $cp=0$, $dp=0$, $a+b=b+a+c$, $a+c=c+a$, $b+c=c+b$, $d+c=c+d$, $a+d=d+a$ and $b+d=d+b$. In particular, each element $x\in R$ is uniquely written in the form $x=ax_1+bx_2+cx_3+dx_4$ with coefficients $0\le x_1<p$, $0\le x_2<p$, $0\le x_3<p$ and $0\le x_4<p$.

Since order of the element $a$ is equal to the exponent of group $G$, i.e. $p$, it follows that by Lemma~\ref{exponent} we can assume that $a$ is an identity of $R$, i.e. $ax=xa=x$ for each $x\in R$. Furthermore, for each $x\in R$ there exist coefficients $\alpha(x)$, $\beta(x)$, $\gamma(x)$, $\varphi(x)$, $\lambda(x)$, $\mu(x)$, $\nu(x)$ and $\psi(x)$ such that $xb=a\alpha(x)+b\beta(x)+c\gamma(x)+d\varphi(x)$ and $xd=a\lambda(x)+b\mu(x)+c\nu(x)+d\psi(x)$. It is clear that they are uniquely defined modulo $p$, so that some mappings $\alpha\colon R\to \mathbb Z_{p}$, ${\beta: R\rightarrow Z_{p}}$, ${\gamma: R\rightarrow Z_p}$, ${\varphi: R\rightarrow Z_p}$, ${\lambda: R\rightarrow Z_p}$, ${\mu: R\rightarrow Z_p}$, ${\nu: R\rightarrow Z_p}$ and ${\psi: R\rightarrow Z_p}$ are determined.

By Corollary~\ref{cor_1}, $L$ is the normal subgroup of order $p^3$ or $p^2$ in $R$. Since $L$ consists the derived subgroup of $R^+$ it follows that the generators $b$ and $c$ we can choose such that $c=-a-b+a+b$. If $|L|=p^3$ then $L=\langle b\rangle + \langle c \rangle+ \langle d \rangle$. Since $R^*=R\setminus L$  it follows ${R^*=\{ax_1+bx_2+cx_3+cx_4\mid x_1\not\equiv 0 \; (\!\!\mod p\;)\}}$ and $x=ax_1+bx_2+cx_3+dx_4$ is invertible if and only if $x_1\not\equiv 0 \; (\!\!\mod p\;)$.

Through this section let $R$ be a local nearring with $|R:L|=p$.

\begin{lemma}\label{G4.4}
If $a$ coincides with identity element of $R$, $x=ax_1+bx_2+cx_3+dx_4$, $y=ay_1+by_2+cy_3+dy_4\in R$, $|R:L|=p$, $xb=a\alpha(x)+b\beta(x)+c\gamma(x)+d\varphi(x)$ and $xd=a\lambda(x)+b\mu(x)+c\nu(x)+d\psi(x)$, then
$$\begin{array}{c}
xy=a(x_1y_1)+b(x_2y_1+\beta(x)y_2+\mu(x)y_4)+c(-x_1x_2\binom{y_1}{2}+\\
+x_3y_1+\gamma(x)y_2+x_1\beta(x)y_3+\nu(x)y_4)+d(x_4y_1+\varphi(x)y_2+\psi(x)y_4).~\rm{(*)}
\end{array}$$
Moreover, for the mappings\\
$\begin{array}{c}
\alpha\colon R\to \mathbb Z_{p}, {\beta: R\rightarrow Z_{p}}, {\gamma: R\rightarrow Z_p}, {\varphi: R\rightarrow Z_p}, {\lambda: R\rightarrow Z_p}, {\mu: R\rightarrow Z_p},\\
{\nu: R\rightarrow Z_p}~\mbox{and}~{\psi: R\rightarrow Z_p}~\mbox{the following statements hold}{\rm :}
\end{array}$
\begin{description}
  \item[\rm{(0)}] $\beta(0)\equiv 0\; (\!\!\mod p)$, ${\gamma(0)\equiv 0\; (\!\!\mod p)}$, ${\varphi(0)\equiv 0\; (\!\!\mod p)}$,\\
   ${\gamma(0)\equiv 0\; (\!\!\mod p)}$, ${\mu(0)\equiv 0\; (\!\!\mod p)}$, ${\nu(0)\equiv 0\; (\!\!\mod p)}$ and\\
   ${\psi(0)\equiv 0\; (\!\!\mod p)}$ if and only if the nearring $R$ is zero-symmetric;
  \item[\rm{(1)}] $\alpha(x)\equiv 0\; (\!\!\mod p)$ and $\lambda(x)\equiv 0\; (\!\!\mod p)${\rm ;}
  \item[\rm{(2)}] if $\beta(x)\equiv 0 \; (\!\!\mod p)$, then $x_1\equiv 0 \; (\!\!\mod p)${\rm ;}
  \item[\rm{(3)}] $\beta(xy)\equiv \beta(x)\beta(y)+\mu(x)\varphi(y)\; (\!\!\mod p\;)${\rm ;}
  \item[\rm{(4)}] $\gamma(xy)\equiv \gamma(x)\beta(y)+x_1\beta(x)\gamma(y)+\nu(x)\varphi(y)\; (\!\!\mod p\;){\rm ;}$
  \item[\rm{(5)}] $\varphi(xy)\equiv \varphi(x)\beta(y)+\psi(x)\varphi(y)\; (\!\!\mod p\;){\rm ;}$
  \item[\rm{(6)}] $\mu(xy)\equiv \beta(x)\mu(y)+\mu(x)\psi(y)\; (\!\!\mod p\;){\rm ;}$
  \item[\rm{(7)}] $\nu(xy)\equiv \gamma(x)\mu(y)+x_1\beta(x)\gamma(y)+\nu(x)\psi(y)\; (\!\!\mod p\;){\rm ;}$
  \item[\rm{(8)}] $\psi(xy)\equiv \varphi(x)\mu(y)+\psi(x)\psi(y)\; (\!\!\mod p\;)$.
\end{description}
\end{lemma}

\begin{proof}
If $|L|=p^3$, then $L=\langle b\rangle + \langle c \rangle + \langle d \rangle$. Since $L$ is the $(R,R)$-subgroup in $R^+$ by statement~1) of Lemma~\ref{prop} it follows that $xb\in L$ and $xd\in L$, hence $a\alpha(x)\in L$ and $a\lambda(x)\in L$ for each $x\in R$. Thus $\alpha(x)\equiv 0 \; (\!\!\mod p)$ and $\lambda(x)\equiv 0 \; (\!\!\mod p)$, so we get statement~(1). Substituting the obtained value of $\alpha(x)\equiv 0 \; (\!\!\mod p)$ and $\lambda(x)\equiv 0 \; (\!\!\mod p)$ in statements~(2)--(4) and (6)--(8) from Lemma~\ref{G4.3}, we obtain statement~(3)--(8) of the lemma and the formula for multiplication $(*)$. Putting $y=c$, we get $xc=c(x_1\beta(x))$. Hence, if $\beta(x)\equiv 0 \; (\!\!\mod p)$, then $xc=0$, and so $x\in L$.  Therefore, $x_1\equiv 0 \; (\!\!\mod p)$, as claimed in statement (2). Indeed, statement (0) repeats the statement (0) of Lemma~\ref{G4.3}.
\end{proof}

Next, we give examples of local nearrings.

\begin{lemma}\label{G4.6}
Let $R$ be a local nearring whose additive group of $R^+$ is isomorphic to $G_4$ and $|R:L|=p$. If $x=ax_1+bx_2+cx_3+dx_4$, $y=ay_1+by_2+cy_3+dy_4\in R$, then the mappings ${\beta: R\rightarrow Z_{p}}$, ${\gamma: R\rightarrow Z_p}$, ${\varphi: R\rightarrow Z_p}$, ${\mu: R\rightarrow Z_p}$, ${\nu: R\rightarrow Z_p}$ and ${\psi: R\rightarrow Z_p}$ from multiplication~$\rm{(*)}$ can be one of the following{\rm :}
\begin{description}
  \item[1)] $\beta(x)={x_1}^i$ and $\psi(x)={x_1}^i$ ($0<i<p$), $\gamma(x)=\varphi(x)=\mu(x)=\nu(x)=0${\rm ;}
  \item[2)] $\beta(x)=1$ and $\psi(x)=1$, $\gamma(x)=\varphi(x)=\mu(x)=\nu(x)=0$.
\end{description}
\end{lemma}

\begin{proof}
It is easy to check that the functions from statements 1) and 2) satisfy conditions 2)--8) of Theorem~\ref{G4.4}.
\end{proof}

As a consequence of Lemma~\ref{G4.6} we have the following result.

\begin{theorem}\label{G4.7}
For each odd prime $p$ there exists a local nearring $R$ whose additive group $R^+$ is isomorphic to $G_4$.
\end{theorem}

\textbf{Example 3.} Let $G\cong C_3 \times ((C_3\times C_3)\rtimes C_3)$. If $x=ax_1+bx_2+cx_3+dx_4$ and $y=ay_1+by_2+cy_3+dy_4\in G$ and $(G,+, \cdot)$ is a local nearring, then by Lemma~\ref{G4.6} $``\cdot"$ can be one of the following multiplications:
\begin{description}
  \item[(1)] $x\cdot y=ax_1y_1+b(x_2y_1+x_1y_2)+c(-x_1x_2\binom{y_1}{2}+x_3y_1+x_1^2y_3)+d(x_4y_1+x_1y_4)${\rm ;}
  \item[(2)] $x\cdot y=ax_1y_1+b(x_2y_1+x_1^2y_2)+c(-x_1x_2\binom{y_1}{2}+x_3y_1+y_3)+d(x_4y_1+x_1^2y_4)${\rm ;}
  \item[(3)] $x\cdot y=ax_1y_1+b(x_2y_1+y_2)+c(-x_1x_2\binom{y_1}{2}+x_3y_1+x_1y_3)+d(x_4y_1+y_4)$.
\end{description}

A computer program verified that for $p=3$ the nearring obtained in Lemma~\ref{G4.6} is indeed a local nearring (see Example~1), is deposited on GitHub:

\verb+https://github.com/raemarina/Examples/blob/main/LNR_81-12.txt+

From the package LocalNR and \cite{zenode_625} we have the following number of all non-isomorphic zero-symmetric local nearrings on $G_4$ of orders $81$ and $625$.

\begin{center}
\begin{tabular}{|c|c|}
\hline
$StructureDescription(R^+)$ & $n(R^+)$\\
 \hline
$C_3 \times ((C_3\times C_3)\rtimes C_3)$                     & \textbf{794} \\
 \hline
 $|R:L|=3$                                                    & 782 \\
 \hline
 $|R:L|=9$                                                    & 12 \\
 \hline
$C_5 \times ((C_5\times C_5)\rtimes C_5)$                     &  \textbf{2090} \\
 \hline
  $|R:L|=5$                                                    & 2078 \\
 \hline
 $|R:L|=25$                                                    & 12 \\
 \hline
\end{tabular}
\end{center}

\subsection{The group $G_5$}

\

Let $G_5$ be additively written group from Theorem~\ref{groups}. Then $G_5=\langle a\rangle +\langle b \rangle +\langle c \rangle$ for some elements $a$, $b$ and $c$ of $R$ satisfying the relations $a{p^2}=b{p}=c{p}=0$, $a+b=b+a(1-p)$, $a+c=c+a$ and $b+c=c+b$.

Recall that the exponent of a finite $p$-group is the maximal order of its elements. The following assertion is easily verified.

\begin{lemma}\label{G5.1}
If $x$ is an element of maximal order in $G_5$, then there exist generators $a$, $b$ and $c$ of this group such that $a=x$ and the relations $a{p^2}=b{p}=c{p}=0$, $-b+a+b=a(1-p)$, $a+c=c+a$, $b+c=c+b$ hold.
\end{lemma}

\begin{lemma}\label{G5.2}
For any natural numbers $k$, $r$, $s$ and $t$ in the group $G_5$ the equalities $ck+bs+ar={ar(1+sp)+bs+ck}$ and $(ar+bs+ck)t=ar(t+s\binom{t}{2}p)+bst+ckt$ hold. \end{lemma}

\begin{proof}
Let $q=1+p$. Since $-b+a+b=a(1-p)$, $a+c=c+a$ and $b+c=c+b$ it follows $c+b+a=aq+b+c$, so $ck+bs+ar=arq^s+bs+ck$ for arbitrary integers $k\ge 0$, $r\ge 0$ and $s\ge 0$. Taking into consideration, that
$$q^s=(1+p)^s\equiv 1+sp \; (\!\mod  p^2)$$
by binomial's formula, giving $ck+bs+ar={ar(1+sp)+bs+ck}$. Next, $(ar+bs+ck)t=ar(1+q^s+\dots+q^{s(t-1)})+bst+ckt$ by induction on $t$. Therefore, $1+q^s+\dots+q^{s(t-1)}\equiv 1+(1+sp)+\dots+(1+s(t-1)p)=t+s\binom{t}{2}p  \; (\!\!\mod  p^2)$, thus $(ar+bs+ckt)t=ar(t+s\binom{t}{2}p)+bst+ckt$.
\end{proof}

\subsection{Nearrings with identity whose\\
additive groups are isomorphic to $G_5$}

\

Let $R$ be a nearring with identity whose additive group $R^+$ is isomorphic to $G_5$. Then $R^{+}=\langle a\rangle +\langle b \rangle +\langle c \rangle$ with elements  $a$, $b$ and $c$, where $a$ coincides with identity element of $R$ and the relations $ap^2=bp=cp=0$, $a+b=b+a(1-p)$, $a+c=c+a$, $b+c=c+b$ are valid.  Moreover, each element $x\in R$ is uniquely written in the form $x=ax_1+bx_2+cx_3$ with coefficients $0\le x_1< p^2$, $0\le x_2<p$ and $0\le x_3<p$.

Consider $a$ coincides with identity element of $R$, so that $xa=ax=x$ for each $x\in R$. Furthermore, for each  $x\in R$ there exist integers $\alpha(x)$, $\beta(x)$, $\gamma(x)$, $\nu(x)$, $\mu(x)$ and $\phi(x)$ such that $xb=a\alpha(x)+b\beta(x)+c\gamma(x)$ and $xc=a\nu(x)+b\mu(x)+c\phi(x)$. It is clear that modulo $p^2$, $p$, $p$ and $p^2$, $p$, $p$, respectively, these integers are uniquely determined by $x$ and so some mappings $\alpha: R\rightarrow\mathbb Z_{p^2}$, $\beta: R\rightarrow\mathbb Z_p$, $\gamma: R\rightarrow\mathbb Z_p$, $\nu: R\rightarrow\mathbb Z_{p^2}$, $\mu: R\rightarrow\mathbb Z_p$ and $\phi: R\rightarrow\mathbb Z_p$ are determined.

\begin{lemma}\label{G5.3} Let $x=ax_1+bx_2+cx_3$ and $y=ay_1+by_2+cy_3$ be elements of $R$. If $a$ coincides with identity element of $R$, then
$$\begin{array}{l}
xy=a(x_1y_1+\alpha(x)y_2+x_1x_2\binom{y_1}{2}p+\nu(x)y_3)+\\
\qquad {}+b(x_2y_1+\beta(x)y_2+\mu(x)y_3)+c(x_3y_1+\gamma(x)y_2+\phi(x)y_3).~\rm{(**)}\end{array}$$
Moreover, for the mappings\\
$$\begin{array}{c}
\alpha: R\rightarrow\mathbb Z_{p^2}, \beta: R\rightarrow\mathbb Z_p, \gamma: R\rightarrow\mathbb Z_p$, $\nu: R\rightarrow\mathbb Z_{p^2},\\
\mu: R\rightarrow\mathbb Z_p~\mbox{and}~\phi: R\rightarrow\mathbb Z_p~\mbox{the following statements hold}{\rm :}
\end{array}$$
\begin{itemize}
  \item[(0)] $\alpha(0)=\beta(0)=\gamma(0)=\nu(0)=\mu(0)=\phi(0)=0$ if and only if the nearring $R$ is zero-symmetric{\rm ;}
  \item[(1)] $\alpha(a)=0$, $\beta(a)=1$, $\gamma(a)=0$, $\nu(c)=0$, $\mu(c)=0$ and $\phi(c)=1${\rm ;}
  \item[(2)] $\alpha(x)\equiv 0\; (\!\!\mod p)$ and $\nu(x)\equiv 0\; (\!\!\mod p)${\rm ;}
  \item[(3)] $\alpha(xy)=x_1\alpha(y)+\alpha(x)\beta(y)+x_1x_2\binom{\alpha(y)}{2}p+\nu(x)\gamma(y)\;(\!\!\mod p^2)${\rm ;}
  \item[(4)] $\beta(xy)=x_2\alpha(y)+\beta(x)\beta(y)+\mu(x)\gamma(y)\;(\!\!\mod p)${\rm ;}
  \item[(5)] $\gamma(xy)=x_3\alpha(y)+\gamma(x)\beta(y)+\phi(x)\gamma(y)\;(\!\!\mod p)${\rm ;}
  \item[(6)] $\nu(xy)=x_1\nu(y)+\alpha(x)\mu(y)+x_1x_2\binom{\mu(y)}{2}p+\nu(x)\phi(y)\;(\!\!\mod p^2)${\rm ;}
  \item[(7)] $\mu(xy)=x_2\nu(y)+\beta(x)\mu(y)+\mu(x)\phi(y)\;(\!\!\mod p)${\rm ;}
  \item[(8)] $\phi(xy)=x_3\nu(y)+\gamma(x)\mu(y)+\phi(x)\phi(y)\;(\!\!\mod p)$.
\end{itemize}\end{lemma}

\begin{proof}
By the left distributive law, we have
$$\begin{array}{c}
xy=(xa)y_1+(xb)y_2+(xc)y_3=(ax_1+bx_2+cx_3)y_1+(a\alpha(x)+b\beta(x)+c\gamma(x))y_2+\\
+(a\nu(x)+b\mu(x)+c\phi(x))y_3.
\end{array}$$
Furthermore, Lemma~\ref{G5.2} implies that $$\begin{array}{l}(ax_1+bx_2+cx_3)y_1=ax_1(y_1+x_2\binom{y_1}{2}p)+bx_2y_1+cx_3y_1,\end{array}$$ $$\begin{array}{l}(a\alpha(x)+b\beta(x)+c\gamma(x))y_2=a\alpha(x)(y_2+\beta(x)\binom{y_2}{2}p)+b\beta(x)y_2+c\gamma(x)y_2,\end{array}$$
and $$\begin{array}{l}(a\nu(x)+b\mu(x)+c\phi(x))y_3=a\nu(x)(y_3+\mu(x)\binom{y_3}{2}p)+b\mu(x)y_3+c\phi(x)y_3.\end{array}$$

By Lemma~\ref{G5.2}, we have
$$\begin{array}{l}bx_2y_1+a\alpha(x)( y_2+\beta(x)\binom{y_2}{2}p)\\
\qquad {}=a\alpha(x)(y_2-\beta(x)\binom{y_2}{2}p)(1-x_2y_1p)+bx_2y_1\end{array}$$

and

$$\begin{array}{l}b\beta(x)y_2+a\nu(x)(y_3+\mu(x)\binom{y_3}{2}p)\\
\qquad {}=a\nu(x)(y_2-\mu(x)\binom{y_2}{2}p)(1-\beta(x)y_2)+b\beta(x)y_2.\end{array}$$

Thus we get
$$\begin{array}{l}xy=a((x_1y_1+\alpha(x)y_2)+(x_1x_2\binom{y_1}{2}+\alpha(x)y_2+\alpha(x)x_2y_1y_2p+\\
\qquad {}+\alpha(x)\beta(x)\binom{y_2}{2})p+\nu(x)y_3+\nu(x)x_2y_1y_2p+\nu(x)\beta(x)y_2y_3p+\\
\qquad {}+\nu(x)\mu(x)\binom{y_3}{2}p)+b(x_2y_1+\beta(x)y_2+\mu(x)y_3)+\\
\qquad {}+c(x_3y_1+\gamma(x)y_2+\phi(x)y_3).\end{array}$$

As $0\cdot a=a\cdot 0=0$, the nearring $R$ is zero-symmetric if and only if $0=0\cdot b=a\alpha(0)+b\beta(0)+c\gamma(0)$ and $0=0\cdot c=a\nu(0)+b\mu(0)+c\phi(0)$ whence $\alpha(0)=\beta(0)=\gamma(0)=\nu(0)=\mu(0)=\phi(0)=0$. Similarly, from the equalities $b=ab=a\alpha(a)+b\beta(a)$ and $c=ac=a\nu(c)+b\mu(c)+c\phi(c)$ it follows that $\alpha(a)=0$, $\beta(a)=1$, $\gamma(a)=0$, $\nu(c)=0$, $\mu(c)=0$ and $\phi(c)=1$, we obtain statement~(1). Since $(xb)p=x(bp)=0$ and $xb=a\alpha(x)+b\beta(x)+c\gamma(x)$, we have $0=(a\alpha(x)+b\beta(x)+c\gamma(x))p=a\alpha(x)(p+\beta(x)\binom{p}{2}p)=a\alpha(x)p$ by Lemma~\ref{G5.2} and hence $\alpha(x)\equiv 0\; {(\!\!\mod p)}$. Moreover, $(xc)p=x(cp)=0$ and $xc=a\nu(x)+b\mu(x)+c\phi(x)$, we have $0=(a\nu(x)+b\mu(x)+c\phi(x))p=a\nu(x)(p+\mu(x)\binom{p}{2}p)=a\nu(x)p$ by Lemma~\ref{G5.2} and hence $\nu(x)\equiv 0\; {(\!\!\mod p)}$, and so statement~(2).
Therefore we obtain
$$\begin{array}{l}xy=a((x_1y_1+\alpha(x)y_2+x_1x_2\binom{y_1}{2}p+\nu(x)y_3)+\\
\qquad {}+b(x_2y_1+\beta(x)y_2+\mu(x)y_3)+c(x_3y_1+\gamma(x)y_2+\phi(x)y_3),\end{array}$$ as desired in~$(**)$.

Finally, the associativity of multiplication in $R$ implies that $x(yb)=(xy)b=a\alpha(xy)+b\beta(xy)+c\gamma(xy)$ and $x(yc)=(xy)c=a\nu(xy)+b\mu(xy)+c\phi(xy)$. Furthermore,  substituting $yb=a\alpha(y)+b\beta(y)+c\gamma(y)$ instead of $y$ in formula~$(**)$, we also have
$$\begin{array}{l}x(yb)=a((x_1\alpha(y)+\alpha(x)\beta(y)+x_1x_2\binom{\alpha(y)}{2}p+\nu(x)\gamma(y))+\\
\qquad {}+b(x_2\alpha(y)+\beta(x)\beta(y)+\mu(x)\gamma(y))+c(x_3\alpha(y)+\gamma(x)\beta(y)+\phi(x)\gamma(y)).\end{array}$$
Comparing the coefficients under $a$ and $b$ in two expressions obtained for $x(yb)$, we derive statements~(3)--(5) of the lemma.

Next, substituting $yc=a\nu(y)+b\mu(y)+c\phi(y)$ instead of $y$ in formula~$(**)$, we get
$$\begin{array}{l}x(yc)=a((x_1\nu(y)+\alpha(x)\mu(y)+x_1x_2\binom{\mu(y)}{2}p+\nu(x)\phi(y))+\\
\qquad {}+b(x_2\nu(y)+\beta(x)\mu(y)+\mu(x)\phi(y))+c(x_3\nu(y)+\gamma(x)\mu(y)+\phi(x)\phi(y)).\end{array}$$
Finally, comparing the coefficients under $a$ and $b$ in two expressions obtained for $x(yc)$, we derive statements~(6)--(8) of the lemma.
\end{proof}

\subsection{Local nearrings whose additive groups\\
are isomorphic to $G_5$}

\

Let $R$ be a local nearring whose additive group $R^+$ is isomorphic to $G_5$. Then $R^{+}=\langle a\rangle +\langle b \rangle +\langle c \rangle$ with elements  $a$, $b$ and $c$, where $a$ coincides with identity element of $R$ and the relations $ap^2=bp=cp=0$, $a+b=b+a(1-p)$, $a+c=c+a$ and $b+c=c+b$ are valid.  Moreover, each element $x\in R$ is uniquely written in the form $x=ax_1+bx_2+cx_3$ with coefficients $0\le x_1< p^2$, $0\le x_2<p$ and $0\le x_3<p$.

Consider $a$ coincides with identity element of $R$, so that $xa=ax=x$ for each $x\in R$. Furthermore, for each  $x\in R$ there exist integers $\alpha(x)$, $\beta(x)$, $\gamma(x)$, $\nu(x)$, $\mu(x)$ and $\phi(x)$ such that $xb=a\alpha(x)+b\beta(x)+c\gamma(x)$ and $xc=a\nu(x)+b\mu(x)+c\phi(x)$. It is clear that modulo $p^2$, $p$, $p$ and $p^2$, $p$, $p$, respectively, these integers are uniquely determined by $x$ and so some mappings $\alpha: R\rightarrow\mathbb Z_{p^2}$, $\beta: R\rightarrow\mathbb Z_p$, $\gamma: R\rightarrow\mathbb Z_p$, $\nu: R\rightarrow\mathbb Z_{p^2}$, $\mu: R\rightarrow\mathbb Z_p$ and $\phi: R\rightarrow\mathbb Z_p$ are determined.

By Corollary~\ref{cor_1}, $L$ is the normal subgroup of order $p^3$ or $p^2$ in $R$. Through this section let $R$ be a local nearring with $|R:L|=p$.

If $|L|=p^3$, then $L=\langle ap\rangle + \langle b \rangle+ \langle c \rangle$. Since $R^*=R\setminus L$  it follows that
$${R^*=\{ax_1+bx_2+cx_3\mid x_1\not\equiv 0 \; (\!\!\mod p\;)\}}$$
and $x=ax_1+bx_2+cx_3$ is invertible if and only if $x_1\not\equiv 0 \; (\!\!\mod p\;)$. Since $L$ is the $(R,R)$-subgroup in $R^+$ by statement~1) of Lemma~\ref{prop} it follows that $xb\in L$ and $xc\in L$, hence $a\alpha(x)\in L$ and $a\nu(x)\in L$ for each $x\in R$. Thus $\alpha(x)\equiv 0 \; (\!\!\mod p)$ and $\nu(x)\equiv 0 \; (\!\!\mod p)$, as in statement~(2) of Theorem \ref{G5.3}. Therefore, for local nearrings $R$ we have the same multiplication as for nearrings with identity, i.e. multiplication $(**)$.

\begin{lemma}\label{G5.3} Let $x=ax_1+bx_2+cx_3$ and $y=ay_1+by_2+cy_3$ be elements of $R$ and $|R:L|=p$. If $a$ coincides with identity element of $R$, then multiplication $(**)$ holds for the mappings from Theorem~\ref{G5.3}.\end{lemma}

Next, we will give examples of local nearrings.

\begin{lemma}\label{G5.4}
Let $R$ be a local nearring whose additive group of $R^+$ is isomorphic to $G_5$ and $|R:L|=p$. If $x=ax_1+bx_2+cx_3$, $y=ay_1+by_2+cy_3\in R$, then the mappings $\alpha: R\rightarrow\mathbb Z_{p^2}$, $\beta: R\rightarrow\mathbb Z_p$, $\gamma: R\rightarrow\mathbb Z_p$, $\nu: R\rightarrow\mathbb Z_{p^2}$, $\mu: R\rightarrow\mathbb Z_p$ and $\phi: R\rightarrow\mathbb Z_p$ can be one of the following{\rm :}
\begin{description}
  \item[1)] $\beta(x)=\phi(x)=\left\{
               \begin{array}{ll}
                 1, & if~\hbox{$x_1\not\equiv 0\; (\!\!\mod p\;);$}\\
                 0, & if~\hbox{$x_1\equiv 0\; (\!\!\mod p\;)$,}
               \end{array}
             \right. \alpha(x)=\gamma(x)=\mu(x)=\nu(x)=0;$
  \item[2)] $\beta(x)=\phi(x)=1$, $\alpha(x)=\gamma(x)=\mu(x)=\nu(x)=0$.
\end{description}
\end{lemma}

\begin{proof}
It is easy to check that the functions from statements 1) and 2) satisfy conditions 1)--8) of Theorem~\ref{G4.4}.
\end{proof}

As a consequence of Lemma~\ref{G5.4} we have the following result.

\begin{theorem}\label{G5.5}
For each odd prime $p$ there exists a local nearring $R$ whose additive group $R^+$ is isomorphic to $G_5$.
\end{theorem}

\textbf{Example 4.} Let $G\cong (C_9\rtimes C_3)\times C_3$. If $x=ax_1+bx_2+cx_3$ and $y=ay_1+by_2+cy_3\in G$ and $(G,+,\cdot)$ is a local nearring, then by Lemma~\ref{G4.6} $``\cdot"$ can be one of the following multiplications.
\begin{description}
  \item[(1)] $x\cdot y=a(x_1y_1+3x_1x_2\binom{y_1}{2})+b(x_2y_1+\beta(x)y_2)+c(x_3y_1+\phi(x)y_3)$, where
  $$\beta(x)=\phi(x)=\left\{
               \begin{array}{ll}
                 1, & if~\hbox{$x_1\not\equiv 0\; (\!\!\mod 3\;)$;}\\
                 0, & if~\hbox{$x_1\equiv 0\; (\!\!\mod 3\;)$,}
               \end{array}\right.$$
  \item[(2)] $x\cdot y=a(x_1y_1+3x_1x_2\binom{y_1}{2})+b(x_2y_1+y_2)+c(x_3y_1+y_3)$.
\end{description}

A computer program verified that for $p=3$ the nearring obtained in Lemma~\ref{G5.4} is indeed a local nearring (see Example~2), is deposited on GitHub:

\verb+https://github.com/raemarina/Examples/blob/main/LNR_81-13.txt+

From the package LocalNR and \cite{zenode_625} we have the following number of all non-isomorphic local nearrings on groups $G_5$ of orders $81$ and $625$.

\begin{center}
\begin{tabular}{|c|c|}
\hline
$StructureDescription(R^+)$ & $n(R^+)$\\
 \hline
$(C_9\times C_3)\rtimes C_3$                        & 337 \\
 \hline
$(C_{25}\times C_5)\rtimes C_5$                     & 630 \\
   \hline
\end{tabular}
\end{center}

\subsection{The groups $G_6$}

\

Let $G_6$ be additively written group from Theorem~\ref{groups}. Then $G_6=\langle a\rangle +\langle b \rangle +\langle c \rangle$ for some elements $a$, $b$ and $c$ of $R$ satisfying the relations $ap^2=0$, $bp=0$, $cp=0$, $a+b=b+a$, $a+c=c+a$, $c+b=b+c+ap$.

\textbf{Conjecture 1.} There does not exist a local nearring whose additive group is isomorphic to $G_6$.

Using GAP and the package LocalNR, Conjecture~1 was confirmed for groups of orders $3^4=81$ and $5^4=625$.

\footnotesize{{\bf Acknowledgement.} The authors are grateful IIE-SRF for their support of our fellowship at the University of Warsaw. This work is partially supported by the Thematic Research Programme ``Tensors: geometry, complexity and quantum entanglement", University of Warsaw, Excellence Initiative --- Research University and the Simons Foundation Award No. 663281 granted to the Institute of Mathematics of the Polish Academy of Sciences for the years 2021--2023.

\end{document}